\documentclass[11pt]{article}
\usepackage{amsmath}
\usepackage{dsfont}
\usepackage{mathrsfs}
\usepackage{amsmath,amssymb}
\usepackage{amsfonts}
\usepackage{hyperref}
\usepackage{amsthm}
\usepackage{graphicx}
\usepackage{subfigure}
\usepackage{xcolor}
\usepackage{cite}
\usepackage{epic,eepic,epsf,epsfig}
\usepackage{floatflt}
\hfuzz=\maxdimen
\theoremstyle{definition}
\newtheorem{lemma}{Lemma}[section]
\newtheorem{definition}[lemma]{Definition}

\newtheorem{theorem}[lemma]{Theorem}

\numberwithin{equation}{section}

\DeclareFixedFont{\Acknowledgment}{OT1}{cmr}{bx}{n}{14pt}
\allowdisplaybreaks

\begin{document}
\title{\textbf{Combinatorial Calabi flows with ideal circle patterns}}
\author{Xiaoxiao Zhang\\  School of Systems Science and Statistics, Beijing Wuzi University, Beijing 100044, P.R. China \\Email: xxzhang0408@126.com}

\date{}
\maketitle

\vspace{-12pt}

\begin{abstract}
In this paper, we extend the work of Ge-Hua-Zhou \cite{GHZ} on combinatorial Ricci flows for ideal circle patterns to combinatorial Calabi flows in both hyperbolic and Euclidean background geometry. We prove that the solution to the combinatorial Calabi flows with any given initial Euclidean (hyperbolic resp.)ideal circle pattern exists for all time and converges exponentially fast to a flat cone metric (hyperbolic resp.) on a given surface.

\end{abstract}
\vspace{12pt}

\section{Introduction}
\quad\quad In modern geometry, researchers strive to find canonical metrics on a given manifold, which are metrics with specific properties that reveal the geometric characteristics of the manifold. To this end, E. Calabi proposed the renowned Calabi flow method\cite{EC1,EC2}, which aims to evolve the manifold until it converges to  constant curvature metric. Specifically, in the study of smooth surfaces, it has been proven that the Calabi flow exists for all time and converges to a constant scalar curvature metric conformal to the initial metric as time tends to infinity(refer to \cite{SC1,SC2,XC,BC,HR}).On the other hand, circle patterns (or circle packings) are another important tool for studying geometry and topology. This method was initially proposed by Koebe\cite{PK} and Andreev\cite{EMA1,EMA2}, who studied arrangements of circles that are tangent to each other in pairs. Thurston(\cite{WT} Chapter 13) further generalized this concept by introducing circle patterns where pairs of circles intersect at acute or right angles, as well as hyperbolic circle patterns for studying the geometry and topology of three-dimensional manifolds. Through hyperbolic circle patterns on triangulated surfaces, Thurston constructed hyperbolic cone metrics with discrete Gaussian curvature, which may have singularities at vertices, thereby providing an in-depth understanding of the geometric characteristics of manifolds.

Chow-Luo, in \cite{CL}, first unveiled a profound connection between Thurston's circle patterns and Hamilton's Ricci flow on surfaces. They introduced the combinatorial Ricci flow on surfaces, a discrete analog of Hamilton's Ricci flow, offering a new perspective in geometric studies. Inspired by Chow-Luo's pioneering efforts, Ge \cite{Ge2,Ge} introduced the combinatorial Calabi flow to deform Thurston's circle patterns, which is exactly the negative gradient flow of the discrete Calabi energy, providing a powerful mathematical tool for understanding the dynamic evolution of circle patterns. Ge's research demonstrated that, in the context of Euclidean circle patterns, the combinatorial Calabi flow exists for all time and converges exponentially fast if and only if Thurston's combinatorial conditions are satisfied. Subsequently, Ge and Xu collaborated in \cite{GX2} to extend this research framework to the combinatorial Calabi flow with hyperbolic circle patterns. They certified that the flow produces solutions which converge to ZCCP-metric (zero curvature circle-packing metric) if the initial energy is small enough, infusing new vitality into hyperbolic geometry research. In \cite{GH}, Ge-Hua collaborated to eliminate the additional assumptions in \cite{GX2}, thus obtaining convergence results for the hyperbolic Calabi flow in generality.

Thurston (\cite{WT}, Chap 13) also observed that ideal circle patterns are closely related to ideal hyperbolic polyhedrons. A convex hyperbolic polyhedron in $\mathbb{B}_3$ is called an idea if all its vertices are on the sphere $\mathbb{S}^2 = \partial\mathbb{B}_3$ at infinity. Given an ideal hyperbolic polyhedron, the boundaries of the hyperbolic planes containing its faces is an ideal circle pattern, its incidence graph is combinatorially dual to the combinatorics of the ideal hyperbolic polyhedron, and the exterior dihedral angles between the faces are equal to the intersection angles between the circles (refer to \cite{BB,JMS,MRo}). Based on this relationship, Ge-Hua-Zhou \cite{GHZ} extended Thurston's characterization about circle pattterns to ideal circle patterns. They showed the uniqueness and the existence of ideal circle patterns, by proving the rigidity (or, say,
injectivity) of the curvature map $K: \mathbb{R}^{|V|}\rightarrow\mathbb{R}^{|V|}$ and describing the image set of $K$. They also derived varies sufficient and necessary conditions, to describe the existence of ideal circle patterns compatible with a flat (hyperbolic resp.) metric on a given surface and used Chow-Luo's surface combinatorial Ricci flows to find the corresponding ideal circle patterns.

Inspired by the main works in \cite{Ge,GH,GHZ}, we consider combinatorial Calabi flows to the ideal circle pattern settings, which is defined as Definition \ref{ideal-c-p}. We will prove the solution to the combinatorial Calabi flow with any given initial Euclidean (hyperbolic resp.) ideal circle pattern exists for all time and converges exponentially fast to a flat cone metric (hyperbolic resp.) on a given surface.

\section{Preliminaries}

Firstly, we recall a natural procedure to convert cellular decomposition $\mathfrak{D}$ into triangulation $\mathcal{T}(\mathfrak{D})$. Given a 2-cell $\tau\in F$, denote its vertices by $v_1, v_2, \ldots, v_m$ in turn. Add a vertex $\tau^{*}$ to $\tau$ and connect it to each $v_i$. Performing this procedure for all 2-cells, we obtain a triangulation $\mathcal{T}(\mathfrak{D})$ of $S$. Let $V^{\diamond}, E^{\diamond}$ and $F^{\diamond}$ be the sets of vertices, edges, and triangles of $\mathcal{T}(\mathfrak{D})$. Then
$$V^{\diamond}=V\cup V^{*},\ E^{\diamond} = E\cup E_{v^{*}v},$$
where $V^{*}$ denotes the vertex set of the dual graph $G^{*}$ of $G$, the set $E_{v^{*}v}$ consists of the edges
$v^{*}v$ such that $v^{*}$ corresponds to a 2-cell containing $v$ as a vertex. Without ambiguity, here and hereafter we do not distinguish the vertex $v\in V$ with its corresponding primal vertex in $V^{\diamond}$.
\subsection{Ideal circle patterns}
\begin{figure}[h]
\centering
\begin{minipage}[t]{0.4\textwidth}
\centering
\includegraphics[width=2in]{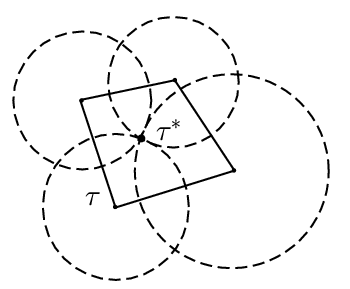}
\caption{circles meet interior}\label{in}
\end{minipage}
\quad\quad\quad
\begin{minipage}[t]{0.35\textwidth}
\centering
\includegraphics[width=1.9in]{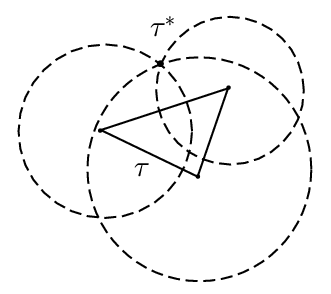}
\caption{circles meet exterior}\label{out}
\end{minipage}
\end{figure}

\quad\quad To better understand the definition of the ideal circle pattern, we recall the definition of the $\mathfrak{D}$-type circle pattern.
Given an oriented closed surface $S$ endowed with a constant curvature metric $\mu$. In this paper, we only concern the Euclidean and the hyperbolic geometry. Assume that $\mathfrak{D}$ is a cellular decomposition of $S$ with sets of vertices, edges and 2-cells $V$, $E$, $F$ respectively. Let $\mathcal{P}=\{C_v,v\in V\}$ be a circle pattern on $(S, \mu)$ with the circles indexed by the vertex set $V$ of a graph $G$, where $G$ is the 1-skeleton of the cellular decomposition $\mathfrak{D}=\{V,E,F\}$ of $S$. Via the geodesic segments homotopic to the edges of $G$, we connect the centers of the corresponding pairs of circles
and obtain a graph $G(\mathcal{P})$, called the \emph{geodesic nerve} of $\mathcal{P}$. The pattern $\mathcal{P}$ is called \emph{$\mathfrak{D}$-type} if $G(\mathcal{P})$ is an embedding graph isotopic to $G$. (Two embedding graph $G_0, G_1$ of $S$ are called \emph{isotopic}, if there exists a continuous map $\Psi: S\times[0,1]\mapsto S$ such that each $\Psi_t=\Psi(\cdot,t)$ is a homeomorphism, $\Psi_0=\text{id\ and}$ $\Psi_1(G_0)=G_1$.) Let $\Theta:E\mapsto[0,\pi)$ be a weight associated the edge set of $G$. A $\mathfrak{D}$-type circle pattern is said realizing the weight function $\Theta$ if the exterior
intersection angle of $C_{v_1}$ and $C_{v_2}$ is equal to $\Theta(e)$, whenever there exists an edge $e\in E$ connecting to the pair of vertices $v_1$ and $v_2$ in $V$.
\begin{definition}\label{ideal-c-p}
A $\mathfrak{D}$-type circle pattern $\mathcal{P}$ on a surface $(S, \mu)$ is called \emph{ideal } if all circles $C_{v_1},C_{v_2},\ldots, C_{v_m}$ meet at an interior\footnote{In Figure\ref{in} (Figure\ref{out} resp.), the circles meet at an interior (exterior resp.) common point.} common point, whenever $v_1,v_2,\ldots, v_m$ form the vertices
of a 2-cell of the cellular decomposition $\mathfrak{D}$.
\end{definition}

 Actually, it is a parallel concept for an ideal circle patterns in \cite{GHZ}. For convenience, we briefly note that an ideal circle pattern is an ideal $\mathfrak{D}$-type circle pattern in this paper. The geometric condition \emph{meet at an interior common point} of ideal is equivalent to\\
 $\langle\star\rangle$\quad For any edge $e, \Theta(e)\in(0,\pi)$, and for any 2-cell in $\mathfrak{D}$ with boundary edges $e_1, e_2,\ldots, e_m$, there holds
$$\sum_{i=1}^m\Theta(e_i)=(m-2)\pi.$$

Obviously, for every $e\in E$, $\Theta(e)\in(0,\pi)$ is necessary true for an ideal circle pattern, and the condition $\langle\star\rangle$ is also necessary. Conversely, assume that $\Theta: E\rightarrow (0,\pi)$ is a function satisfying $\langle\star\rangle$. When $S$ is of genus $g=0$, Rivin \cite{IR} proved the unique existence of an ideal $\mathfrak{D}$-type circle pattern. When $S$ is of genus $g>0$, Bobenko-Springborn \cite{BAS} certified the unique existence of an ideal $\mathfrak{D}$-type circle pattern. It means that there exists a unique point $\tau^{*}$ in the interior of any $\tau=v_1v_2\cdots v_m\in F$, such that the geodesic distance from $v_i$ to $\tau^{*}$ is exactly $r_i$ for all $1\leq i\leq m$, i.e., the circles $C_{v_1} ,\ldots , C_{v_m}$ centered at $v_1, \ldots , v_m$ with radius $r_1,\ldots , r_m$ intersect at $\tau^{*}$. So they form an ideal circle pattern. We may consider the unique common intersection point $\tau^{*}$
as the geometric dual of the 2-cell $\tau$.

\section{Statement of results}

 \quad\quad Given an oriented closed surface $S$ with a cellular decomposition $\mathfrak{D}=\{V,E,F\}$, where $V$, $E$, $F$ represent the set of vertices, edges and 2-cells respectively. To get an ideal circle pattern, one may consider an alternative way. Let $\Theta: E\rightarrow (0,\pi)$ be an edge weight satisfying the condition $\langle\star\rangle$. For any $r: V \rightarrow(0,+\infty)$, we endow $S$ a piecewise flat (hyperbolic resp.) cone metric as follows. Endow each edge $\{ij\}\in E$ a length
\begin{equation}\label{cos-law-H}
  l_{ij}=\cosh^{-1}(\cosh r_i \cosh r_j + \sinh r_i\sinh r_j\cos\Theta_{ij})
\end{equation}
in the hyperbolic geometry background and
\begin{equation}\label{cos-law-E}
l_{ij}=\sqrt{r_i^2+r_j^2+2r_ir_j\cos\Theta_{ij}}
\end{equation}
in the Euclidean geometry background. For each 2-cell $\tau\in F$, there exists a Euclidean (hyperbolic resp.) polygon with edge length between any two adjacent vertices $i$ and $j$ is exactly $l_{ij}$ determined by (\ref{cos-law-E}) ((\ref{cos-law-H}) resp.). To show this, one first corresponds a Euclidean
(hyperbolic resp.) triangle with three lengths $r_i, r_j$ and $l_{ij}$ to each boundary edge $\{ij\}$ of $\tau$. The existence of these triangles is guaranteed by Lemma \ref{p4}. One then get such a polygon by gluing all these triangles coherently along the common edges, so that these triangles
have a common vertex $\tau^{*}$. Note the geodesic distance from $\tau^{*}$ to a vertex $i\in\partial\tau$ is exactly $r_i$, and the condition $\langle\star\rangle$ is essentially used here to guarantee the existence of $\tau^{*}$. Gluing all these Euclidean (hyperbolic resp.) polygon along the common edges, we obtain a hyperbolic (or Euclidean) cone metric $\mu(r)$ on $S$ with possible cone singularities
at vertices of $T(\mathfrak{D})$.

It is noticeable that at each star-vertex $\tau^{*}$, which is the dual of a 2-cell in $\mathfrak{D}$ with boundary edges $e_1, e_2,\ldots, e_m$, the cone angle is $$\sum_{i=1}^m(\pi-\Theta(e_i))=m\pi-(\sum_{i=1}^m\Theta(e_i))=m\pi-(m-2)\pi=2\pi.$$
Therefore, the cone singularity at any added vertex is naturally smooth. It remains to consider the cone singularities at the primal vertices. For each $v_i \in V$, define the corresponding discrete curvature by
\begin{equation*}
 K_i=2\pi-\text{the\ cone\ angle\ at}\ v_i.
\end{equation*}

For $i,j\in V$, we often write $i\thicksim j$ if the vertices $i$ and $j$ are adjacent. Any function defined on $V$ can be regarded as a column vector in $\mathbb{R}^{|V|}$, hence we consider a circle pattern $r:=(r_1,r_2,\ldots,r_{|V|})^T$ as a point in $\mathbb{R}_{>0}^{|V|}$. A circle pattern $r$ determines the curvature function $K=(K_1,\ldots,K_{|V|})^T$. Conversely, $K$ also determines $r$. In fact, from Andreev-Thurston rigidity on hyperbolic circle patterns (refer to \cite{EMA1,EMA2,dV}), we know a hyperbolic circle pattern is uniquely determined by its curvatures, i.e., the curvature map $K(r): \mathbb{R}_{>0}^{|V|}\rightarrow \mathbb{R}^{|V|}, r\mapsto K(r)$ is injective. For any ideal Euclidean (or hyperbolic) circle pattern, Ge-Hua-Zhou \cite{GHZ} gives analogous results that is the curvature map $K$ is injective.

\subsection{Hyperbolic geometry background}

\quad\quad Now, we recall the definition of combinatorial Calabi flow in \cite{GH}.
Given an oriented closed surface $S$ with a triangulation $\mathcal{T(\mathfrak{D})}$. For each ideal circle pattern $r\in\mathbb{R}_{>0}^{|V|}$, the corresponding \emph{discrete Calabi energy} is
\begin{equation}
\mathcal{C}(r)=\|K\|^2=\sum_{i=1}^{|V|}K_i^2.
\end{equation}
Set $u_i=\ln\tanh\frac{r_i}{2}$, for any $i\in V$. Then the map $u=u(r)$ is a diffeomorphism between $\mathbb{R}_{>0}^{|V|}$ and $\mathbb{R}_{<0}^{|V|}$. The gradient (with respect to $u$-coordinate) of the discrete Calabi energy is
\begin{equation}
\nabla_u\mathcal{C}=(\nabla_{u_1}\mathcal{C},\ldots,\nabla_{u_{|V|}}\mathcal{C})^T=2L K,
\end{equation}
where $K=(K_1,K_2,\ldots,K_{|V|})^T$ and $L$ is the Jacobian of the curvature map $K=K(u)$. The \emph{combinatorial Calabi flow} is defined as
\begin{equation}
 u'(t)=-LK=-\frac{1}{2}\nabla_u\mathcal{C}.
\end{equation}

\begin{theorem}\label{exists-for-all-time}
Let $G$ be the 1-skeleton of a cellular decomposition $\mathfrak{D}$ of an oriented closed surface $S$. Assume the edge weight $\Theta: E\rightarrow (0,\pi)$ satisfies $\langle\star\rangle$. Then the solution $r(t)$ to the Calabi flow with any given initial circle pattern $r(0)$ exists for all time $t\geq0$ and converges exponentially fast to a hyperbolic metric on $S$ with no singularity if and only if there exists an ideal circle pattern with zero curvature.
\end{theorem}

\subsection{Euclidean geometry background}
Firstly, we recall the two-circle configurations on Euclidean in \cite{GHZ} . Let $C_i, C_j$ be two Euclidean circles with vertex $i, j$ and radius $r_i, r_j$, respectively. Assume $C_i$ and $C_j$ intersect at a point $\tau^{*}$. Parallelling to Lemma \ref{p4}, we have a configuration of two intersection circles in $\mathbb{E}^2,$ unique up to Euclidean isometry(i.e. congruence), having radii $\{r_i,r_j\}$ and intersection angles $\Theta_{ij}$. In this Euclidean triangle $\Delta ij\tau^{*}$, denote $\theta_i$ ($\theta_j$ resp.) as the inner angle at the vertex $i$. ($j$ resp.)

\begin{theorem}\label{Euclidean}
Let $G$ be the 1-skeleton of a cellular decomposition $\mathfrak{D}$ of an oriented closed surface $S$. Assume the edge weight $\Theta: E\rightarrow (0,\pi)$ satisfies $\langle\star\rangle$. Consider the combinatorial Calabi flow under Euclidean geometry background. For any given initial pattern $r(0)$, the solution to the combinatorial Calabi flow uniquely exists for all the time $t\geq 0$. Additionally, $r(t)$ converges exponentially fast to some constant curvature circle pattern metric if and only if there exists the constant curvature ideal circle pattern metric.
\end{theorem}

The rigidity and existence of ideal circle patterns has already been resolved by Bobenko-Springborn \cite{BS} and Ge-Hua-Zhou \cite{GHZ}. For the convenience of readers, we append Ge-Hua-Zhou's related work as below.

\begin{theorem}(Ge-Hua-Zhou).
Assume that  $\Theta: E\rightarrow (0,\pi)$ satisfies the condition $\langle\star\rangle$. In hyperbolic background geometry, the solution $r(t)$ to the Ricci flow $\frac{dr_i}{dt}=-K_i\sinh r_i$ for $i=1,2,\cdots,|V|$, exists for all time $t\geq 0$, and the following $H_1-H_5$ are all equivalent:\\
$H_1.$ $r(t)$ converges as $t\rightarrow+\infty$ .\\
$H_2.$ The origin $(0,0,\cdots,0)$ belongs to the image of the curvature map.\\
$H_3.$ When $A$ is a non-empty subset of $V$, the following condition holds:
$$\sum_{e\in E,\partial e\cap A\neq\emptyset}\Theta(e)>\pi |A|.$$
$H_4.$ When $A$ is a non-empty subset of $V$, the following condition holds:
$$\sum_{(e,\tau)\in LK(A)}(\Theta(e)-\pi)+2\pi\chi^2(S(A))<0.$$
$H_5.$ The genus $g>1$ and the weight satisfies $\sum_{l=1}^s\Theta(e_l)<(s-2)\pi$, whenever $e_1,e_2,\cdots,e_n$ form a pseudo-Jordan curve which is not the boundary of a 2-cell of $\mathfrak{D}$.
\end{theorem}

\begin{theorem}(Ge-Hua-Zhou).
The Euclidean background geometry, the solution $r(t)$ to the Ricci flow $\frac{dr_i}{dt}=(K_{av}-K_i)r_i$ for $i=1,2,\cdots,|V|$, exists for all time $t\geq 0$, and the following $E_1-E_4$ are all equivalent:\\
$E_1.$ $r(t)$ converges as $t\rightarrow+\infty$ .\\
$E_2.$ The vector $(K_{av},K_{av},\cdots,K_{av})$ belongs to the image of the curvature map.\\
$E_3.$ When $A$ is a proper non-empty subset of $V$, the following condition holds:
$$\frac{\pi\chi^2(S(A))}{|V|}>\pi |A|-\sum_{e\in E,\partial e\cap A\neq\emptyset}\Theta(e).$$
$E_4.$ When $A$ is a proper non-empty subset of $V$, the following condition holds:
$$\frac{2\pi\chi^2(S(A))}{|V|}>\sum_{(e,\tau)\in LK(A)}(\Theta(e)-\pi)+2\pi\chi^2(S(A)).$$
Moreover, if the flow converges, then it converges exponentially fast to a radius vector so that all vertex curvatures are equal to $K_{av}$.
\end{theorem}

\section{Proof of Theorem \ref{exists-for-all-time}}

\quad\quad We know that a two-circle configuration may be considered as the limiting case of three-circle configuration as one radius of three circles tends to 0. Hence an ideal circle pattern may be considered as the limiting case of the circle pattern. Let $C_i$ and $C_j$ be two hyperbolic circles at vertex $i$ and $j$, and with radius $r_i$ and $r_j$ respectively. Assume $C_i$ and $C_j$ intersects at a point $\tau^{*}$ with an angle $\Theta_{ij}\in(0,\pi)$. \\

Now, we introduce some lemmas which are useful to our main results.

\begin{figure}[h]
\centering
\begin{minipage}[t]{0.4\textwidth}
\centering
\includegraphics[width=2in]{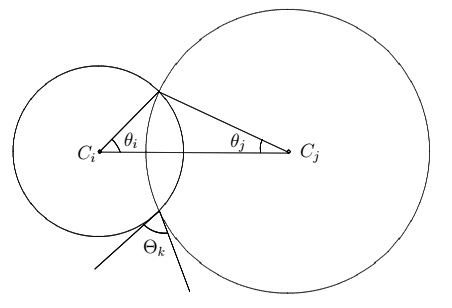}
\caption{pattern of two circles}\label{two-circle}
\end{minipage}
\quad\quad\quad
\begin{minipage}[t]{0.35\textwidth}
\centering
\includegraphics[width=2.5in]{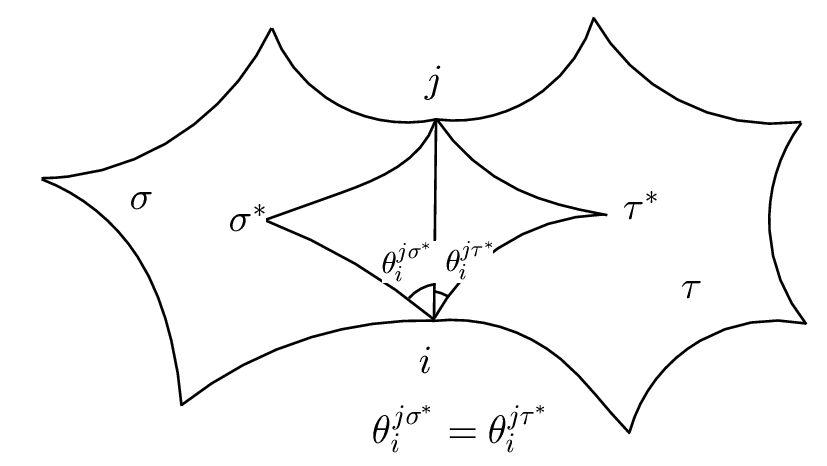}
\caption{two adjacent cells}\label{two-adjacent-cells}
\end{minipage}
\end{figure}

\begin{lemma}\label{p4}
Given $\Theta_k\in(0,\pi)$, any two positive numbers $r_i,r_j$, there exists a configuration of two intersection circles in both Euclidean and hyperbolic geometries, unique up to isometry, having radii $\{r_i,r_j\}$ and meeting in exterior intersection angles $\Theta_k$ (see Figure \ref{two-circle}).
\end{lemma}

Now assume the edge weight $\Theta:E\mapsto(0,\pi)$ satisfies $\langle\star\rangle$. For each edge ${ij}\in E$, let $\sigma$
and $\tau$ be the two adjacent 2-cells with $\{ij\}$ as a common edge, see Figure \ref{two-adjacent-cells}. Let ${\sigma}^{*}$ ($\tau^{*}$ resp.) be the geometric dual of $\sigma$ ($\tau$ resp.), then we see $l_{i\sigma^{*}}=l_{i\tau^{*}}=r_i,\ l_{j\sigma^{*}}=l_{j\tau^{*}}=r_j$ and $\theta_{\sigma^{*}}^{ij}=\theta_{\tau^{*}}^{ij}=\Theta_{ij}$. Hence the hyperbolic triangles $\Delta ij\sigma^{*}$ and $\Delta ij\tau^{*}$ are congruent. Note the angles $\theta_i^{j\sigma^{*}}$ and $\theta_i^{j\tau^{*}}$ are equal, and they depend only on $r_i, r_j$, and $\Theta_{ij}$. We may well denote these two angles as
$$\theta_i^j=\theta_i^{j\sigma^{*}}\ (=\theta_i^{j\tau^{*}})$$
and denote the area of the hyperbolic triangle $\triangle ij\sigma^{*}$ as
$$ \text{Area}_{ij}=\text{Area}(\triangle ij\sigma^{*})(= \text{Area}(\Delta ij\tau^{*}))$$
Thus we can express the curvature at $i$ as
\begin{equation}\label{curvature}
  K_i=2\pi-2\sum_{j\thicksim i}\theta_i^j.
\end{equation}

The following lemma is parallel to those results in Thurston (\cite{WT},
Chap 13), Marden-Rodin \cite{MR}, Colin de Verdi$\grave{e}$re \cite{dV}, Chow-Luo \cite{CL}, Guo-Luo \cite{GL}, Guo
\cite{RG}, Zhou \cite{ZZ} and others, respectively.

\begin{lemma} (\cite{GHZ}) Let $\triangle_{ij\tau^{*}}$ be a triangle which is patterned by two circles with fixed weight $\Theta_{ij}\in (0,\pi)$ as intersection angles. Let $\theta_i$ be the inner angle at vertex $i$. Then
\\ (i) in Euclidean geometry background,
\begin{equation}\label{p7}
  \frac{\partial\theta_i}{\partial r_i}<0;\quad \frac{\partial \theta_i}{\partial r_j}>0;\quad \frac{\partial \theta_i}{\partial r_j}r_j=\frac{\partial \theta_j}{\partial r_i} r_i.,
\end{equation}
\\ (ii) in hyperbolic geometry background,
\begin{equation}\label{p1}
  \frac{\partial\theta_i}{\partial r_i}<0;\quad \frac{\partial \theta_i}{\partial r_j}>0;\quad \frac{\partial \text{Area} (\triangle_{ij\tau^{*}})}{\partial r_i}>0,
\end{equation}
\begin{equation}\label{p2}
\frac{\partial \theta_i}{\partial r_j}\sinh r_j=\frac{\partial \theta_j}{\partial r_i}\sinh r_i.
\end{equation}
\end{lemma}

The following useful result is similar to the related results in Ge-Jiang \cite{GJ} and Ge-Xu \cite{GX}.
\begin{lemma}(\cite{GHZ})\label{p3}
In the hyperbolic triangle $\triangle ij\tau^{*}$, for any $\epsilon>0$, there exists a number $l$ such that whenever $r_i>l$, the inner angle $\theta_i$ is smaller than $\epsilon$.
\end{lemma}
The result in Lemma \ref{p3} is similar to the related results in Ge-Jiang \cite{GJ} and Ge-Xu \cite{GX}.

\begin{lemma}\label{Area-e}
There exists a large enough constant number $C>0$ such that $r_i>C$, then
$$\frac{\partial (\text{Area}_{ij}-\theta_j)}{\partial r_i}\sinh r_i>0.$$
\end{lemma}
\noindent\begin{proof}
Set $c_{ij}=\min(\cos\Theta_{ij},0)$, then $-1<c_{ij}\leq 0$ and $\cosh r_i\cosh r_j>e^{r_i+r_j-\ln 4}.$ Hence we have
\begin{align*}
e^{l_{ij}}\geq\cosh l_{ij}&=\cosh r_{i}\cosh r_{j}+\sinh r_{i}\sinh r_{j}\cos \Theta_{ij}\\&=(1+\cos \Theta_{ij})\cosh r_{i}\cosh r_{j}-\cos \Theta_{ij}\cosh(r_{i}-r_{j})\\&\geq(1+c_{ij})\cosh r_{i}\cosh r_{j}\\&\geq(1+c_{ij})e^{r_i+r_j-\ln 4}.
\end{align*}
That is
\begin{equation*}
  l_{ij}\geq r_i+r_j-\ln\frac{4}{1+c_{ij}}.
\end{equation*}
In view of $r_i>C$ and $C$ is a large enough number, we get
\begin{equation}\label{cosh-bound}
  \cosh l_{ij}>2.
\end{equation}
By the hyperbolic sine law, we have
$$\sinh l_{ij}\sin \theta_i=\sinh r_{j}\sin \Theta_{ij}.$$
Differential the above equality with respect to $r_i$, we get
$$\frac{\partial \theta_i}{\partial r_i}\sinh r_i=-\frac{\partial \theta_i}{\partial r_j}\sinh r_j\cosh l_{ij}.$$
Using the area formula in hyperbolic geometry, we get the area of the hyperbolic triangle $\triangle ij\tau^{*}$ is $\text{Area}(\triangle ij\tau^{*})=\pi-(\theta_i+\theta_j+\pi-\Theta_{ij})=\Theta_{ij}-\theta_i-\theta_j$.
\begin{align}
\frac{\partial (\text{Area}_{ij}-\theta_j)}{\partial r_i}\sinh r_i&=-\frac{\partial \theta_i}{\partial r_i}\sinh r_i-2\frac{\partial \theta_j}{\partial r_i}\sinh r_i\nonumber\\&=\frac{\partial \theta_j}{\partial r_i}\sinh r_i(\cosh l_{ij}-2)
\end{align}
By lemma \ref{p1}, we know $\frac{\partial \theta_j}{\partial r_i}>0$. Combining this proposition with (\ref{cosh-bound}), there holds
$$\frac{\partial (\text{Area}_{ij}-\theta_j)}{\partial r_i}\sinh r_i>0.\quad\quad\Box$$
\end{proof}

Under the condition $\langle\star\rangle$, for each vertex $i$, its curvature $K_i$ depends only on $r_i$, all radius $r_j$ and all weights $\Theta_{ij}$ with $j\thicksim i$. Using Lemma \ref{p2}, we have
$$\frac{\partial K_i}{\partial u_j}=-2\frac{\partial \theta_i^j}{\partial u_j}<0,$$ $$ \frac{\partial K_i}{\partial u_i}=-\sum_{j\thicksim i}\frac{\partial K_i}{\partial u_j}+2\sum_{j\thicksim i}\frac{\partial \text{Area}_{ij}}{\partial u_i}.$$
Set $A=diag\{A_1,A_2,\ldots,A_{|V|}\}$, $L_B=((L_B)_{ij})_{|V|\times |V|}$, where $A_i=2\sum_{j\thicksim i}\frac{\partial Area_{ij}}{\partial r_i}\sinh r_i$ and
\begin{align*}
 (L_B)_{ij}=\begin{cases}-\sum_{k\thicksim i}\frac{K_i}{\partial u_k},\quad &\text{if}\ j=i;\\ \quad \frac{\partial K_i}{\partial u_j},&\text{if}\ j\thicksim i;\\ \quad\quad 0,&\text{if}\ j\not\thicksim i\ \text{and}\ j\neq i.\end{cases}
\end{align*}
Then the Jacobian matrix $L=\partial(K_1,K_2,\ldots,K_{|V|})/\partial(u_1,u_2,\ldots,u_{|V|})$ of $K(u)$ satisfies $L=A+L_B$.

\begin{lemma}(\cite{GHZ})\label{p5}
The Jacobian matrix of $K(u)$ (as a function of $u$) is positive definite.
\end{lemma}

\begin{theorem}\label{bouned-above}
Let $G$ be the 1-skeleton of a cellular decomposition $\mathfrak{D}$ of an oriented closed surface $S$. Assume the edge weight $\Theta: E\rightarrow (0,\pi)$ satisfies $\langle\star\rangle$. Let $r(t)$ be the unique solution to the ideal Calabi flow on a maximal the interval $[0,T)$. Then all $r_i(t)$ are uniformly bounded above on $[0,T)$.
\end{theorem}
\noindent\begin{proof}
By contradiction. Suppose it is not true, then there exists at least one vertex $i\in V$, such that
\begin{equation}
\limsup_{t\rightarrow T}r_i(t)=+\infty.
\end{equation}
For this vertex $i$, using Lemma \ref{p3}, we can choose a large enough positive number $l$ such that $r_i>l$, the inner angle $\theta_i$ is smaller than $\frac{\pi}{2d_i}$, where $d_i$ is the degree of the vertex $i$. Then we have $K_i>\pi$.\\
Set $L=\max\{l,c,r_i(0)+1\}$, where $c$ is given in Lemma \ref{Area-e}. Now, we claim that for any $t\in(0,T)$ and if $r_i(t)>l$, then
\begin{equation}\label{<0}
\frac{dr_i}{dt}<0.
\end{equation}
Since
\begin{align*}
  \sum_{j\thicksim i}\frac{\partial K_i}{\partial u_j}+A_i&=-2\sum_{j\thicksim i}\frac{\partial \theta_i}{\partial u_j}+2\sum_{j\thicksim i}\frac{\partial \text{Area}_{ij}}{\partial r_i}\sinh r_i\\&=-2\sum_{j\thicksim i}\frac{\partial \theta_j}{\partial r_i}\sinh r_i+2\sum_{j\thicksim i}\frac{\partial \text{Area}_{ij}}{\partial r_i}\sinh r_i
  \\&=2\sum_{j\thicksim i}\frac{(\partial \text{Area}_{ij}-\theta_j)}{\partial r_i}\sinh r_i\\&>0.
\end{align*}
we have
\begin{align*}
\frac{1}{\sinh r_i}\frac{d r_i}{d t}&=\Delta K_i=-(L_B+A)K_i\\
&=\sum_{j\thicksim i}\frac{\partial K_i}{\partial u_j}(K_i-K_j)-A_iK_i
\\&<-\pi \sum_{j\thicksim i}\frac{\partial K_i}{\partial u_j}-A_iK_i
\\&<-\pi(\sum_{j\thicksim i}\frac{\partial K_i}{\partial u_j}+A_i)\\&<0.
\end{align*}
Hence
\begin{equation}\label{dr<0}
  \frac{d r_i}{d t}<0.
\end{equation}
By (\ref{<0}), we may choose $t_0\in(0,T)$ such that $r_i(t_0)>c$. Let $t_1\in[0,t_0]$ attain the maximum of $r_i(t)$ in $[0,t_0]$. By the definition of $L$, $t_1>0$. Hence
$$\frac{d r_i}{d t}(t_1)\geq 0,$$
which contradicts to (\ref{dr<0}). This proves the theorem.\quad\quad$\Box$
\end{proof}

\

\noindent\textbf{Proof of Theorem \ref{exists-for-all-time}}
$\Rightarrow$ Since the solution $r(t)$ to the ideal Calabi flow converges, that is, $r(t)$ converges in the Euclidean space topology to some $r^{*}\in\mathbb{R}_{>0}^{|V|}$ as time $t$ tends to $+\infty$. Let $u^{*}=\ln\tanh \frac{r^{*}}{2}$, then $u(t)$ converges to $u^{*}$.
$$u(n+1)-u(n)=u'(\xi_n)=-LK(\xi_n)\rightarrow 0$$
as $n\rightarrow+\infty$, where $\xi_n$ is some number in $(n,n+1)$. As $r(t)\rightarrow r^{*}$, we have $K(r)\rightarrow K(r^{*})$. Thus $LK(r^{*})=0$. Since $L$ is a positive definite matrix, we have $K(r^{*})=0$. The uniqueness of the ideal circle pattern with zero curvature from that $K$ is injective.

$\Leftarrow$ Since $K$ is injective, we let $\bar{r}$ be the unique ideal circle pattern with $K(\bar{r})=0$. Let $r(t)$ be the unique solution to the ideal Calabi flow on a maximal time interval $[0,T)$, we need to prove $T=+\infty$ and $r(t)\rightarrow \bar{r}$ exponentially fast.\\
Let $\bar{u}\in \mathbb{R}_{<0}^{|V|}$ be the $u$-coordinate of $\bar{r}$. Consider the ideal Ricci potential
\begin{equation}\label{line-int}
F(u)\triangleq\int_{\bar{u}}^u\sum_{i=1}^{|V|}K_idu_i,\ u\in\mathbb{R}_{>0}^{|V|}.
\end{equation}
This type of line integral was first introduced by de Verdi\`{e}re \cite{dV}. By the fact $\partial K_i/\partial u_j=\partial K_j/\partial u_i$, the smooth differential 1-form $\sum_{i=1}^{|V|}K_idu_i$ is closed. Hence the line integral (\ref{line-int}) is well defined and is independent on the choice of piecewise smooth paths in $\mathbb{R}_{<0}^{|V|}$ from $u_0$ to $u$.
\\ By lemma 6.1 in [?], there holds
\begin{equation}\label{infinite}
\lim_{\|u\|\rightarrow+\infty, u\in\mathbb{R}_{<0}^{|V|}}F(u)=+\infty.
\end{equation}
By lemma \ref{p5}, that is, the Jacobian matrix $L$ of $K(u)$ is positive definite, we have
\begin{equation*}
\frac{d}{dt}F(u(t))=-\sum_{i}K_iL K_i=-K^TLK\leq 0,
\end{equation*}
which implies that $F(u(t))$ is descending and bounded from below. Hence $F(u(+\infty))$ exists. Moreover, there exists a time sequence $\xi_n\rightarrow+\infty$ such that
\begin{equation}\label{000}
F(u(n+1))-F(u(n))=\frac{d}{dt}\mid_{t=\xi_n}F(u(t))=-K^TLK(\xi_n)\rightarrow 0
\end{equation}
as $n$ goes to $+\infty$.
By (\ref{infinite}), there exists a positive constant $M$ such that $u_i(t)\geq-M$ for all $i$ and $t$ and $M$ is only dependent on the triangulation $\mathcal{T}$ and the initial circle pattern $r(0)$. It follows that $r_i(t)\geq C>0$ for all $i$ and $t$, where $C=\ln\frac{1+e^{-M}}{1-e^{-M}}$ is a constant depending only on the triangulation $\mathcal{T}$ and the initial ideal circle pattern $r(0)$. Hence, by Theorem \ref{bouned-above}, we see $\{r(t)\mid t\in[0,\infty)\}$ is compactly supported in $\mathbb{R}_{>0}^{|V|}.$ Thus there exists an ideal circle pattern $\hat{r}\in\mathbb{R}_{>0}^{|V|}$, and a subsequence $t_k=\xi_{n_k},k\geq 1$ such that $r(t_k)\rightarrow \hat{r}$ as $k$ goes to $+\infty$. Because $K=K(r)$ is a continuous map of $r$, so $K(r(t_k))\rightarrow \hat{r}$ as $k$ goes to $+\infty$. From the limit (\ref{000}) and positive definite of the matrix $L$, we see $K(\hat{r})=0$. By the rigidity result,
i.e., the curvature map $K$ is injective, $\hat{r}=\bar{r}$ is the unique ideal circle pattern with zero curvature.

Since $\{r(t)\mid t\in[0,\infty)\}$ is compactly supported in $\mathbb{R}_{>0}^{|V|}$, we have $\{u(t)\mid t\in[0,\infty)\}\subset\subset\mathbb{R}_{<0}^{|V|}$. Then the eigenvalue of matrix $L$ has a uniform positive lower bound $\lambda_1$, that is $\lambda(L(t))\geq\lambda_1$ along the ideal Calabi flow. In view of the fact $\frac{dK}{dt}=-L^2K$, we have
$$\mathcal{C}'(t)=2K^T\dot{K}(t)=-2K^TL^2K\leq-2\lambda_1^2K^TK=-2\lambda_1^2\mathcal{C}(t).$$
So $\mathcal{C}(t)\leq\mathcal{C}(0)e^{-2\lambda_1^2t}$ and $|K_i(t)|\leq|K(t)|=\sqrt{\mathcal{C}(t)}\leq\sqrt{\mathcal{C}(0)}e^{-\lambda_1^2t}$.
As $\{u(t)\mid t\in[0,\infty)\}\subset\subset\mathbb{R}_{<0}^{|V|}$, we know $L_{ij}$ and $\sinh r_i$ are bounded along the ideal Calabi flow. Hence
\begin{align*}
|\frac{du_i}{dt}|&=|\Delta K_i|=|\sum_{j}L_{ij}K_j|\leq c_1e^{-\lambda_1^2t},\\
|\frac{dr_i}{dt}|&=|\sinh r_i\Delta K_i|=|\sinh r_i\sum_{j}L_{ij}K_j|\leq c_2e^{-\lambda_1^2t},
\end{align*}
where $c_1, c_2$ are positive constants. This implies that the solution converges
with exponential rate.\quad\quad$\Box$

\section{Proof of Theorem \ref{Euclidean}}
Let $C_i$ and $C_j$ be two Euclidean circles with vertex $i$ and $j$, and radius $r_i$ and $r_j$ respectively. Assume $C_i$ and $C_j$ intersect at a point $\tau^{*}$ with an angle $\Theta_{ij}\in(0,\pi)$. Parallelling to Lemma \ref{p4}, we have a configuration of two intersecting circles in $\mathbb{E}^2$, unique up to Euclidean isometry (i.e. congruence), having radii $\{r_i,r_j\}$ and intersection angles $\Theta_{ij}$. In this Euclidean triangle $\triangle ij\tau^{*}$, denote $\theta_i$ ($\theta_j$ resp.) as the inner angle at the vertex $i$ ($j$ resp.). By the Euclidean cosine law, we have $l_{ij}=\sqrt{r_i^2+r_j^2+2r_ir_j\cos\Theta_{ij}}$. Parallelling to Lemma \ref{p1},
Lemma \ref{p2}, we have
\begin{equation}
  \frac{\partial \theta_i}{\partial r_j}r_j=\frac{\partial \theta_j}{\partial r_i}r_i=\frac{d_{ij}}{l_{ij}}>0,
\end{equation}
where $d_{ij}$ is the distance from $\tau^{*}$ to the geodesic line passing $i$ and $j$. Moreover,
\begin{equation}
  \frac{\partial \theta_i}{\partial r_i}r_i=\frac{\partial \theta_j}{\partial r_j}r_j=-\frac{d_{ij}}{l_{ij}}<0.
\end{equation}
Set $\omega_{ij}=d_{ij}/{l_{ij}} > 0$ for each edge $i\thicksim j$. Using the coordinate change $u_i = \ln r_i$ for each
vertex $i\in V$.

\noindent\textbf{Proof of Theorem \ref{Euclidean}}
 Let $d_i$ denotes the degree at vertex $i\in V$. Suppose $d$ be the maximum degree at all vertices, i.e., $d=\max\{d_1,d_2,\ldots,d_{|V|}\}$. Then for any $i\in V$, we have $(2-2d)\pi<K_i<2\pi$.
\\ Now we claim all $|\omega_{ij}|$ are uniformly bounded by a constant $c(\Theta)$ which is only dependent on $\Theta$.
Since
\begin{align*}
  l_{ij}&=r_i^2+r_j^2+2r_ir_j\cos\Theta_{ij} \\&=r_i^2+2r_ir_j\cos\Theta_{ij}+r_j^2\cos^2\Theta_{ij}+r_j^2\sin^2\Theta_{ij}
  \\ &\geq r_j^2\sin^2\Theta_{ij},
\end{align*}
we have $l_{ij}\geq r_j\sin\Theta_{ij}$. Then
$$\omega_{ij}=\frac{d_{ij}}{l_{ij}}\leq\frac{r_j}{l_{ij}}\leq\frac{1}{\sin\Theta_{ij}}.$$
Because it is a finite decomposition on a oriented closed surface, $\sin\Theta_{ij}$ has a positive lower bound $c(\Theta)$. That is to say $|\omega_{ij}|\leq c(\Theta)$.
So
$$|\frac{du_i(t)}{dt}|\leq|\Delta K_i|\leq\sum_{j\thicksim i}|w_{ij}||K_j-K_i|\leq 2d\pi\sum_{j\thicksim i}|w_{ij}|\leq 2d^2\pi c(\Theta)=c_1,$$
we have $-c_1t\leq u_i(t)\leq c_1t$. Hence,
$$c_0e^{-c_1t}\leq r_i(t)\leq c_0e^{c_1t},$$
where $c_0=c(r(0))$, Which implies the ideal combinatorial Calabi flow has a solution for all time $y\in[0,\infty)$ for any $r(0)\in \mathbb{R}_{>0}^{|V|}.$
\quad\quad$\Box$

\noindent{\bf Acknowledgments:}
 The author is supported by National Natural Science Foundation of China under Grant No. 12301069. I would like to thank Professor H. Ge for telling me that Shengyu Li and Zhigang Wang at Hunan First Normal university have also independently obtained the similar results.

\end{document}